\documentclass{ifacconf}

\counterwithin*{section}{part}

\usepackage{amssymb,mathtools}
\usepackage{natbib}
\usepackage{cleveref}
\usepackage{tikz}

\newcommand{\X}{\mathcal X}
\newcommand{\U}{\mathcal U}
\newcommand{\Y}{\mathcal Y}
\newcommand{\V}{\mathcal V}
\newcommand{\K}{\mathcal K}
\newcommand{\M}{\mathcal M}
\newcommand{\N}{\mathcal N}
\newcommand{\Z}{\mathcal Z}
\newcommand{\W}{\mathcal W}
\newcommand{\Rf}{\mathbf R}
\newcommand{\adjun}{^{\ast}}
\newcommand{\Flow}{\mathcal F}
\newcommand{\Effort}{\mathcal E}
\newcommand{\D}{\mathcal D}
\newcommand{\Lagr}{\mathcal L}
\newcommand{\Rc}{\mathcal R}
\DeclareMathOperator{\dom}{dom}
\DeclareMathOperator{\graph}{graph}
\DeclareMathOperator{\Rea}{Re}
\DeclarePairedDelimiterX{\dualprod}[2]{\langle}{\rangle}{#1,#2}
\DeclarePairedDelimiterX{\scprod}[2]{\langle}{\rangle}{#1,#2}
\DeclarePairedDelimiterX{\setdef}[2]{\{}{\}}{#1\,\delimsize\vert\,\mathopen{} #2}
\newcommand{\pvek}[2]{\begin{pmatrix}#1\\#2\end{pmatrix}}
\newcommand{\spvek}[2]{\left(\begin{smallmatrix}#1\\#2\end{smallmatrix}\right)}
\newcommand{\pre}{\operatorname{pr}_{\mathrm e}}
\newcommand{\prf}{\operatorname{pr}_{\mathrm f}}
\crefname{thm}{Theorem}{Theorems}
\crefname{lem}{Lemma}{Lemmas}
\crefname{prop}{Proposition}{Propositions}
\crefname{cor}{Corollary}{Corollaries}
\crefname{figure}{Fig.}{Figs.}
\crefname{table}{Table}{Tables}

\begin{document}
\begin{frontmatter}
\title{Algebraic Constraint Transfers in Infinite-Dimensional Port-Hamiltonian Systems\thanksref{funding}}
\thanks[funding]{Funded by the Deutsche Forschungsgemeinschaft (DFG,
German Research Foundation), Project-ID 531152215, CRC 1701.}
\author[ilmenau]{Anna Fischer}
\author[ilmenau]{Timo Reis}
\address[ilmenau]{Institut f\"ur Mathematik, Technische Universit\"at Ilmenau, \\
Weimarer Stra\ss e 25, 98693 Ilmenau, Germany\\
e-mail: \{anna.fischer, timo.reis\}@tu-ilmenau.de}
\begin{abstract}
Algebraic constraints in geometric port-Hamiltonian systems may occur
in the Dirac, Lagrangian, or resistive relation. We transfer selected
Dirac constraints on storage efforts to the Lagrangian relation, Lagrangian state
constraints to the Dirac structure, and resistive effort or flow
constraints to auxiliary Dirac ports. The constructions apply to
infinite-dimensional Hilbert spaces, preserve the relevant adjoint symmetries and power properties, and retain the restrictions imposed by
the domains of unbounded operators. We characterize the correspondence of classical
trajectories and the required temporal regularity. A string with a massless
segment and a heat conductor containing a perfectly conducting region
illustrate the transfers. 
\end{abstract}
\end{frontmatter}

\section{Introduction}
Geometric port-Hamiltonian systems describe energy storage, power-conserving interconnection, and dissipation by means of Lagrangian, Dirac, and resistive relations, respectively.
Algebraic constraints may occur in any of them. We investigate how such constraints can be transferred
between the relations while preserving the power balance,
and how the trajectories of the resulting representations
correspond.

In finite dimensions, \citet{vanderSchaftMaschke2018}
introduced a generalized port-Hamiltonian framework in which
energy storage is described by Lagrangian subspaces. They also
showed how algebraic constraints can be transferred between
Dirac and Lagrangian structures in both directions.
The extension to nonlinear systems was developed in
\citet{vanderSchaftMaschke2020}.

Even for linear systems, these transfers require care in infinite
dimensions: Dirac structures may involve unbounded differential
operators whose domains encode spatial regularity and boundary
conditions. Such domain restrictions differ from algebraic
constraints, represented by closed, possibly infinite-dimensional
subspaces. Replacing an operator domain by its closure may
change the system. Early work on port-Hamiltonian PDE systems with boundary ports
includes \citet{vanderSchaftMaschke2002,LeGorrecZwartMaschke2005,
Villegas2007}.
Further developments address infinite-dimensional Dirac structures
\citep{IftimeRomanSandovici2014}, analytical aspects of
port-Hamiltonian PDEs \citep{BrugnoliHaineMatignon2023},
differential storage relations \citep{MaschkevanderSchaft2020Boundary},
geometric PDE representations
\citep{BendimeradHohlMatignonHaineLefevre2024,
BendimeradHohlHaineLefevreMatignon2025}, boundary conditions
\citep{GernandtPhilippPreusterSchaller2026}, and
structure-preserving discretization
\citep{BendimeradHohlHaineLefevreMatignon2026}.

We transfer selected closed algebraic constraints between Dirac
and Lagrangian structures on possibly infinite-dimensional
Hilbert spaces, allowing unbounded operators. We also transfer
effort and flow constraints of self-adjoint, nonpositive resistive
relations to auxiliary lossless Dirac ports. These constructions
preserve adjoint symmetries, power properties, and operator-domain
restrictions. We characterize the correspondence of classical trajectories, including the necessary regularity conditions.

Converting representations of complete port-Hamiltonian
systems is a related but distinct problem. In finite dimensions,
\citet{CaballeriaRamirezLefevreLeGorrec2026} developed
equivalent state- and effort-based formulations with resistive
and external ports. In \citet{FischerReis2026Descriptor},
we converted geometric and descriptor representations on
Hilbert spaces, expressing Dirac constraints through auxiliary
equations. Here we remain within the geometric framework
and transfer selected constraints between Dirac, Lagrangian,
and resistive relations.

We illustrate the results by two PDE models with boundary
ports: an elastic string containing a massless segment and
a heat conductor with a perfectly conducting region.

\section{Geometric systems and closed constraints}
\subsection{Ports, relations, and trajectories}
All Hilbert spaces are complex, and their inner products are linear in
the first argument. For a Hilbert space $\U$, let $\U\adjun$ be the space of
bounded antilinear functionals on $\U$,
$\Rf_\U:\U\to\U\adjun$ the Riesz isomorphism, and
$\iota_\U:\U\to\U^{\ast\ast}$ the canonical biduality map.
We identify Hilbert spaces with their biduals via this map
whenever convenient. We use the
canonical pairing $\dualprod{u'}{u}_{\U\adjun,\U}$; the pairing with reversed arguments is defined by complex conjugation. 
For a bounded operator
$A:\U\to\V$, its anti-dual adjoint
$A\adjun:\V\adjun\to\U\adjun$ is characterized by
$\dualprod{A\adjun v'}{u}=\dualprod{v'}{Au}$. For mutually orthogonal subspaces $\mathcal A$ and $\mathcal B$,
we write $\mathcal A\perp\mathcal B$ for their algebraic direct sum
$\mathcal A+\mathcal B$.

A linear relation is a subspace of a product space, not necessarily an operator graph, allowing for algebraic constraints and multivalued constitutive laws.

For an operator $A$, we write
\[\graph(A)=\setdef*{(x,Ax)}{x\in\dom A}.\]
For Hilbert spaces $\U,\V$ and a unitary identification
$\Gamma:\V\to \U\adjun$, define
$\dualprod{v}{u}_{\Gamma}=\dualprod{\Gamma v}{u}_{\U\adjun,\U}$.
A linear relation $\mathcal T\subseteq \U\times \V$ has domain
$\dom\mathcal T$ and multivalued part
$\setdef*{v}{(0,v)\in\mathcal T}$. Its adjoint with respect to $\Gamma$ is
\[
 \mathcal T^{*_{\Gamma}}:=\setdef*{(a,c)}{
 \begin{aligned}
  &\dualprod{c}{u}_{\Gamma}
   =\overline{\dualprod{v}{a}_{\Gamma}}\\[-2pt]
  &\text{for all }(u,v)\in\mathcal T
 \end{aligned}}.
\]
The relation $\mathcal T$ is called self-adjoint if
$\mathcal T^{*_{\Gamma}}=\mathcal T$ and skew-adjoint if
\[
 \mathcal T^{*_{\Gamma}}
 =\setdef*{(u,-v)}{(u,v)\in\mathcal T}.
\]
Neither condition requires the relation to be single-valued
or its domain to be dense.

Let $\X$ be the state space and $\Flow_{\mathrm R}$  the
resistive flow space, and put $\Effort_{\mathrm R}=\Flow_{\mathrm R}\adjun$.
The external input and output belong to $\U$ and $\U\adjun$,
respectively. The total flow and effort spaces are
\[
 \Flow_{\mathrm{tot}}=\X\times\Flow_{\mathrm R}\times\U\adjun,
 \qquad
 \Effort_{\mathrm{tot}}=\X\adjun\times\Effort_{\mathrm R}\times\U,
\]
and we write
\begin{multline*}
 \dualprod{e}{f}_{\mathrm{tot}}
 :=\dualprod{e_{\mathrm S}}{f_{\mathrm S}}
   +\dualprod{e_{\mathrm R}}{f_{\mathrm R}}
   +\dualprod{u}{f_{\mathrm P}},\\
   f=\left(\begin{matrix}f_{\mathrm S}\\f_{\mathrm R}\\f_{\mathrm P}\end{matrix}\right),\quad e=\left(\begin{matrix}e_{\mathrm S}\\e_{\mathrm R}\\u\end{matrix}\right).
\end{multline*}
The last term uses the pairing with reversed arguments.

A \emph{Dirac structure}
$\D\subseteq\Flow_{\mathrm{tot}}\times\Effort_{\mathrm{tot}}$
describes the power-conserving interconnection. It is a closed
linear relation equal to its orthogonal complement with respect
to the symmetrized power form
\[
 [(f,e),(g,h)]_+
 :=\dualprod{e}{g}_{\mathrm{tot}}
 +\overline{\dualprod{h}{f}_{\mathrm{tot}}}.
\]
Equivalently, $\D$ is skew-adjoint with respect to
the pairing $\dualprod{\cdot}{\cdot}_{\mathrm{tot}}$.
In particular, $\D$ is maximally isotropic and every
Dirac pair has zero real power.

A \emph{Lagrangian structure}
$\Lagr\subseteq\X\times\X\adjun$ describes admissible
state--storage-effort pairs and is self-adjoint with respect to the pairing
$\dualprod{\cdot}{\cdot}_{\X\adjun,\X}$.
In particular, $\Lagr$ is maximally symmetric.

A resistive relation
$\Rc\subseteq\Flow_{\mathrm R}\times\Effort_{\mathrm R}$
describes dissipative interactions. It is \emph{nonpositive} if
$\Rea\dualprod{e_{\mathrm R}}{f_{\mathrm R}}\leq0$
for all $(f_{\mathrm R},e_{\mathrm R})\in\Rc$. 
 
The complete geometric system is the inclusion
\begin{equation}\label{eq:full-system}
 \begin{gathered}
 \left(
 \begin{pmatrix}\phantom{-}\dot x(t)\\-f_{\mathrm R}(t)\\-y(t)\end{pmatrix},
 \begin{pmatrix}e_{\mathrm S}(t)\\e_{\mathrm R}(t)\\u(t)\end{pmatrix}
 \right)\in\D,\\
 (x(t),e_{\mathrm S}(t))\in\Lagr,\qquad
 (f_{\mathrm R}(t),e_{\mathrm R}(t))\in\Rc.
 \end{gathered}
\end{equation}
A classical trajectory on an interval $I$ has
$x\in C^1(I;\X)$, while all efforts, internal resistive flows,
and external port variables are continuous in their respective spaces;
\eqref{eq:full-system} holds pointwise. The storage functional is
\[
 \mathcal H(x,e_{\mathrm S})=
 \tfrac12\dualprod{e_{\mathrm S}}{x}_{\X\adjun,\X}.
\]
By self-adjointness of $\Lagr$, this functional is real-valued
and independent of the choice of effort for a given admissible
state.
Writing $\mathcal H(t):=\mathcal H(x(t),e_{\mathrm S}(t))$,
the Lagrangian property gives
$\mathcal H(t)-\mathcal H(s)=\tfrac12\Rea\dualprod{e_{\mathrm S}(t)
+e_{\mathrm S}(s)}{x(t)-x(s)}$ for $s,t\in I$. Thus
$\frac{\mathrm d}{\mathrm dt}\mathcal H(t)
=\Rea\dualprod{e_{\mathrm S}(t)}{\dot x(t)}$
without requiring a time derivative of the effort. Dirac isotropy
then yields
\begin{equation}\label{eq:power}
 \tfrac{\mathrm d}{\mathrm dt}\mathcal H(x(t),e_{\mathrm S}(t))
 =\Rea\dualprod{e_{\mathrm R}(t)}{f_{\mathrm R}(t)}+\Rea\dualprod{u(t)}{y(t)}.
\end{equation}
No positivity of storage is needed for the transfer theorems.
Internal resistive and external ports can be omitted
by taking the corresponding spaces to be zero.

The transformations differ in where they impose an additional
zero condition; \cref{tab:transfers} summarizes the four cases.
The auxiliary spaces need not be finite-dimensional.
\begin{table}[!ht]
\caption{Location of the transferred constraint.}
\label{tab:transfers}
\centering
\setlength{\tabcolsep}{4pt}
\begin{tabular}{@{}l l@{\ }r@{}l l@{\ }l l@{\ }l@{}}
\hline
Original & \multicolumn{3}{l}{Zero condition} & \multicolumn{2}{l}{Variable} & \multicolumn{2}{l@{}}{Result} \\
\hline
Dirac effort     & storage state & $z$     & ${}=0$ & effort & $\mu$     & Thm. & \ref{thm:DL} \\
Lagrangian state & Dirac effort  & $e_\mu$ & ${}=0$ & state  & $\mu$     & Thm. & \ref{thm:LD} \\
Resistive effort & Dirac port    & $qe$    & ${}=0$ & effort & $\mu$     & Thm. & \ref{thm:Rtransfer} \\
Resistive flow   & Dirac port    & $qf$    & ${}=0$ & flow   & $\lambda$ & Cor. & \ref{cor:Rflow} \\
\hline
\end{tabular}
\end{table}

\subsection{Relaxation of selected algebraic constraints}
\label{subsec:relation-preparation}
Let $\mathcal T\subseteq\U\times\V$ be self-adjoint
or skew-adjoint with respect to $\Gamma$. We call a closed subspace $\K\subseteq\U$ satisfying
$\dom\mathcal T\subseteq\K$ a \emph{selected algebraic
constraint}. 
For a subspace $\mathcal A\subseteq\U$, we write
$\mathcal A^\perp$ for its orthogonal complement and
\[
 \mathcal A^\circ:=\setdef*{v\in\V}{
 \dualprod{v}{a}_\Gamma=0
 \text{ for all }a\in\mathcal A}
\]
for its annihilator.
Annihilators of subspaces of $\V$ are defined analogously.
By adjoint symmetry, the multivalued part of
$\mathcal T$ is $(\dom\mathcal T)^\circ$.

For a closed subspace $\W$, let $P_{\W}$ denote
the orthogonal projection onto $\W$.
For a flow--effort relation $\mathcal S$, we write
$\prf\mathcal S$ and $\pre\mathcal S$ for its flow
and effort projections, respectively.

\begin{lem}[Relaxing a selected constraint]\label[lem]{lem:relaxation}
Let $\mathcal T\subseteq \U\times \V$ be self-adjoint or skew-adjoint
with respect to $\Gamma$. For a closed subspace $\K\subseteq\U$ containing
$\dom\mathcal T$, set
\[
 \M=\K^\perp,\quad \N=\M^\circ,\quad
 \mathcal T_0=\mathcal T\cap(\U\times \N).
\]
Then
\begin{equation}\label{eq:relax}
 \mathcal T=\mathcal T_0\perp(\{0\}\times \K^\circ),\quad
\widetilde{\mathcal T}\phantom{\vcentcolon}
   \mathrel{\mathllap{\vcentcolon}}=\mathcal T_0\perp(\M\times\{0\})
\end{equation}
are direct sums. The relation $\widetilde{\mathcal T}$ retains
self-adjointness or skew-adjointness; nonpositivity is also
preserved in the self-adjoint case.
Moreover,
\begin{equation}\label{eq:relax-domain}
 \dom\widetilde{\mathcal T}=\dom\mathcal T\perp \M.
\end{equation}
If $\K=\overline{\dom\mathcal T}$, the relaxed relation is the
graph of a densely defined closed operator. If, in addition,
$\K=\dom\mathcal T$, this operator is bounded.

An analogous construction applies to constraints on the
second component. If a closed subspace $\K\subseteq\V$
contains the second projection of $\mathcal T$, set
\[
 \M=\K^\perp,\quad \N=\M^\circ,\quad
 \mathcal T_1=\mathcal T\cap(\N\times\V).
\]
Then
\begin{equation}\label{eq:relax-second}
 \mathcal T=\mathcal T_1\perp(\N^\perp\times\{0\}),\qquad
\widetilde{\mathcal T}\phantom{\vcentcolon}
   \mathrel{\mathllap{\vcentcolon}}=\mathcal T_1\perp(\{0\}\times\M)
\end{equation}
are direct sums. The relation
$\widetilde{\mathcal T}$ retains self-adjointness or
skew-adjointness; nonpositivity is also preserved in
the self-adjoint case.
\end{lem}
\begin{pf}
Since $\dom\mathcal T\subseteq\K$, its multivalued part
contains $\K^\circ$. The unitary $W=\Rf_\U^{-1}\Gamma$
satisfies $W(\K^\circ)=\K^\perp$ and $W(\N)=\K$.
Transporting the second component by $W$ decomposes
$\mathcal T$ into a closed part in $\K\times\K$ and
$\{0\}\times\K^\perp$. Taking adjoints along
$\U=\K\perp\K^\perp$ shows that the first part retains
the adjoint symmetry. The horizontal factor
$\K^\perp\times\{0\}$ is both self-adjoint and skew-adjoint;
replacing the vertical factor by it preserves that symmetry
and nonpositivity, since it contributes zero power.
Equation \eqref{eq:relax-domain} follows; density and the
closed graph theorem give the remaining claims.

For the second-component case, the projection condition
and adjoint symmetry imply
$\N^\perp\times\{0\}\subseteq\mathcal T$, yielding
\eqref{eq:relax-second}. With respect to the induced
duality between $\N^\perp$ and $\M$, this
horizontal factor and $\{0\}\times\M$ are both self-adjoint
and skew-adjoint and contribute zero power. Exchanging
them preserves the adjoint symmetry and nonpositivity.
\end{pf}

For the first-component construction, set
$\Y=\U/\K$ and let $q:\U\to\Y$ be the quotient map.
The bounded map
\begin{equation}\label{eq:reaction}
 j=\Gamma^{-1}q\adjun:\Y\adjun\to \V,
 \qquad
 \dualprod{j\mu}{u}_{\Gamma}
 =\dualprod{\mu}{qu}
\end{equation}
is an isometric isomorphism onto $\K^\circ$, giving the
bounded recovery map
$r_\K=j^{-1}P_{\K^\circ}:\V\to\Y\adjun$.
Each $(u,v)\in\mathcal T$ has the unique decomposition
$v=v_0+j\mu$, where $(u,v_0)\in\mathcal T_0$ and
$\mu=r_\K v$. Both $\K$ and $\Y$ may be infinite-dimensional.
In the second-component construction, the corresponding
decomposition is $u=u_1+P_{\N^\perp}u$, where
$(u_1,v)\in\mathcal T_1$.

\begin{cor}[Partial transfers]
Set $\K_0=\overline{\dom\mathcal T}$ in the
first-component construction of \cref{lem:relaxation}. Then
\[
 \overline{\dom\widetilde{\mathcal T}}
 =\K_0\perp \K^\perp.
\]
The transferred constraints correspond to $\K^\perp$.
The remaining multivalued part is
$(\K_0\perp \K^\perp)^\circ$, corresponding to
$\K\cap \K_0^\perp$. The possibly nonclosed operator
domain within $\K_0$ is retained.
\end{cor}
\begin{pf}
Use \eqref{eq:relax-domain} and the orthogonal decomposition
\[\U=\K_0\perp(\K\cap \K_0^\perp)\perp \K^\perp.\]
\end{pf}

\noindent\emph{A proper dense domain is not an algebraic constraint.}
Consider $A=\frac{\mathrm d}{\mathrm d\xi}$ on $L^2(0,1)$ with
\[\dom A=\setdef*{u\in H^1(0,1)}{u(0)=u(1)}.\] Integration by parts shows that $A$ is skew-adjoint, while
$\dom A\subsetneq L^2(0,1)$ and
$\overline{\dom A}=L^2(0,1)$.
The associated Dirac structure
$\setdef*{(Au,\Rf_{L^2}u)}{u\in\dom A}$ thus has no nontrivial
closed algebraic effort constraint. Its periodic trace condition
belongs to the operator domain and is not removed by a transfer.

\section{Transfers between Dirac and Lagrangian structures}
In both transfers, the resistive and external channels of
\eqref{eq:full-system} remain unchanged. 

\begin{thm}[Dirac to Lagrangian]\label{thm:DL}
Assume that all storage efforts in $\D$ belong to a closed
$\K\subseteq\X\adjun$. Apply the second-component relaxation
of \cref{lem:relaxation} directly to $\D$, using the closed
subspace $\K\times\Effort_{\mathrm R}\times\U$ of the total
effort space. Denote the resulting Dirac structure by
$\widetilde\D$. Put $\Y=\X\adjun/\K$, let $q:\X\adjun\to\Y$ be the quotient map,
and set $j=\iota_{\X}^{-1}q\adjun:\Y\adjun\to\X$ and
$j_{\mathrm S}\mu=\left(\begin{smallmatrix}j\mu\\0\\0\end{smallmatrix}\right)$. Define, up to a permutation of channels,
\begin{align}
 \widehat\D={}&\setdef*{\left(
 \pvek{f+j_{\mathrm S}\mu}{-qe_{\mathrm S}},\pvek e\mu\right)}
 {(f,e)\in\widetilde\D,\ \mu\in\Y\adjun},
 \label{eq:DL-D}\\
 \widehat\Lagr={}&\setdef*{\left(
 \pvek x0,\pvek{e_{\mathrm S}}\mu\right)}
 {(x,e_{\mathrm S})\in\Lagr,\ \mu\in\Y\adjun}.
 \label{eq:DL-L}
\end{align}
Then $\widehat\D$ is Dirac and $\widehat\Lagr$ is Lagrangian.
The selected effort constraint becomes $z=0$, where
$z\in\Y$ is the auxiliary storage state. Classical trajectories correspond
bijectively, preserving $x,f_{\mathrm R},e_{\mathrm R},u,y$.
The auxiliary effort is uniquely determined by
$j\mu=P_{\K^\circ}\dot x$, where
$\K^\circ\subseteq\X$ is the annihilator of $\K$.
\end{thm}
\begin{pf}
The relaxed relation $\widetilde\D$ is Dirac by
\cref{lem:relaxation}. The auxiliary multiplier is related to
the Dirac-constraint representation used in
\citet{FischerReis2026Descriptor}; here the additional zero
condition is placed in the Lagrangian relation. The bounded coupling
$\spvek{e_{\mathrm S}}{\mu}\mapsto\spvek{j\mu}{-qe_{\mathrm S}}$
is skew-adjoint by \eqref{eq:reaction}. The corresponding flow shear maps the Dirac structure
$\widetilde\D\times(\{0\}_{\Y}\times\Y\adjun)$
onto \eqref{eq:DL-D} while preserving the symmetrized
power form. Hence $\widehat\D$ is Dirac. The relation in \eqref{eq:DL-L}
is the product of $\Lagr$ with the vertical Lagrangian
factor $\{0\}_{\Y}\times\Y\adjun$.
Along an extended trajectory, $z=0$ and
$\dot z=-qe_{\mathrm S}$ imply $qe_{\mathrm S}=0$.
The decomposition in \eqref{eq:relax-second} then recovers
the original Dirac relation. Conversely, the same
decomposition uniquely determines the auxiliary multiplier
through $j\mu=P_{\K^\circ}\dot x$. Since $j$ is an isometric
isomorphism onto $\K^\circ$, the multiplier is continuous
whenever $\dot x$ is continuous. All remaining ports
are untouched.
\end{pf}

\begin{thm}[Lagrangian to Dirac]\label{thm:LD}
Let $\K\subseteq\X$ be closed and contain the state projection
of $\Lagr$. Relax $\Lagr$ by \cref{lem:relaxation}, and denote
the resulting relation by $\widetilde\Lagr$. Put
$\Y=\X/\K$, $q:\X\to\Y$, and
$j=q\adjun:\Y\adjun\to\X\adjun$. Define
\begin{align}
 \widehat\Lagr={}&\setdef*{\left(
 \pvek x\mu,\pvek{e+j\mu}{qx}\right)}
 {(x,e)\in\widetilde\Lagr,\ \mu\in\Y\adjun},
 \label{eq:LD-L}\\
 \widehat\D={}&\D\times(\Y\adjun\times\{0\}_{\Y}),
 \label{eq:LD-D}
\end{align}
with the auxiliary storage channel grouped with $\X$.
The relations $\widehat\Lagr$ and $\widehat\D$ are
Lagrangian and Dirac, respectively. The selected state constraint becomes $e_\mu=0$,
where $e_\mu=qx$ is the auxiliary effort. Every classical extended trajectory
projects to a classical original one. Conversely, an original
classical trajectory admits a classical extension if and only if
its uniquely determined multiplier $\mu=r_\K e_{\mathrm S}$ belongs
to $C^1(I;\Y\adjun)$.
\end{thm}
\begin{pf}
The bounded map $\spvek{x}{\mu}\mapsto\spvek{j\mu}{qx}$ is
self-adjoint by \eqref{eq:reaction}. Adding it to the
self-adjoint relation
$\widetilde\Lagr\times(\Y\adjun\times\{0\}_{\Y})$
shows that \eqref{eq:LD-L} is Lagrangian. 
The horizontal auxiliary factor in \eqref{eq:LD-D} is a Dirac structure.
Its vanishing effort
gives $qx=0$, which recovers $\Lagr$ via \eqref{eq:relax}.
The multiplier is uniquely determined by the original effort
and depends continuously on it. Since it is now part of the
state, however, an extended classical trajectory requires
$\dot\mu$ to exist and be continuous.
\end{pf}

\noindent\emph{The regularity condition is essential.}
Take $\X=\mathbb C$ and the (vertical) relations
$\D=\Lagr=\{0\}\times\mathbb C$. The functions
$x(t)=0$, $e_{\mathrm{S}}(t)=|t|$ form an original classical trajectory
on $(-1,1)$. After the Lagrangian-to-Dirac transfer with
$\K=\{0\}$, the auxiliary state is necessarily $\mu(t)=|t|$,
which is not $C^1$. Thus the original trajectory has no
classical extension, although the pointwise relation correspondence
remains exact.

\section{Resistive constraints and lossless auxiliary ports}
We transfer constraints of a self-adjoint, nonpositive
resistive relation
$\Rc\subseteq\Flow_{\mathrm R}\times\Effort_{\mathrm R}$.
Unlike \cref{thm:DL,thm:LD}, these transfers modify
$\Rc$ while preserving both properties and leaving
storage and external ports unchanged.

A bounded skew-adjoint map from efforts to flows defines
the power-preserving \emph{flow shear}
$(F,E)\mapsto(F+BE,E)$; a bounded skew-adjoint map
from flows to efforts gives the \emph{effort shear}
$(F,E)\mapsto(F,E+BF)$. Both are boundedly invertible and preserve the symmetrized
power form, hence they map Dirac structures to Dirac structures. Auxiliary resistive flows enter the Dirac relation
with the same sign convention as in \eqref{eq:full-system}.

\begin{thm}[Resistive effort constraints]\label{thm:Rtransfer}
Let $\K\subseteq\Effort_{\mathrm R}$ be closed and contain
$\pre\Rc$. With $\M=\K^\perp$, $\N=\M^\circ$,
$\Y=\Effort_{\mathrm R}/\K$, $q:\Effort_{\mathrm R}\to\Y$,
and $j=\iota_{\Flow_{\mathrm R}}^{-1}q\adjun$, put
\begin{equation}\label{eq:Rrelax}
 \widetilde\Rc=
 (\Rc\cap(\N\times\Effort_{\mathrm R}))
 \perp(\{0\}\times \M).
\end{equation}
Introduce the bounded lossless coupling
\begin{equation}\label{eq:RB}
 \begin{gathered}
 B:\Effort_{\mathrm R}\times\Y\adjun
 \to\Flow_{\mathrm R}\times\Y,\\
 B\spvek{e}{\mu}=\spvek{j\mu}{-qe}.
 \end{gathered}
\end{equation}
Let $\mathcal B$ extend $B$ by zero to the unchanged
channels, and define
\begin{align}
 \D^{\mathrm a}&=\D\times(\{0\}_{\Y}\times\Y\adjun),\nonumber\\
 \widehat\Rc&=\widetilde\Rc\times
                 (\{0\}_{\Y}\times\Y\adjun),
 \label{eq:Rhat}\\
 \widehat\D&=\setdef*{(F+\mathcal B E,E)}{(F,E)\in\D^{\mathrm a}}.
 \label{eq:RDhat}
\end{align}
Then $\widehat\D$ is Dirac and $\widehat\Rc$ is self-adjoint
and nonpositive. Classical trajectories correspond bijectively, preserving
the storage and external port variables; the resistive
variables satisfy
\begin{equation}\label{eq:Rrecover}
 e_{\mathrm R}=e,\qquad
 f_{\mathrm R}=f+j\mu,\qquad
 qe=0,\quad(f,e)\in\widetilde\Rc.
\end{equation}
No additional dynamic state is introduced. The auxiliary
effort is uniquely determined by
$j\mu=P_{\N^\perp}f_{\mathrm R}$ and is therefore
continuous for classical trajectories.
\end{thm}
\begin{pf}
By the second-component case of \cref{lem:relaxation},
$\widetilde\Rc$ in \eqref{eq:Rrelax} is self-adjoint and
nonpositive. The coupling $B$ in \eqref{eq:RB} is
skew-adjoint by \eqref{eq:reaction}; 
thus the bounded flow shear in \eqref{eq:RDhat}
maps $\D^{\mathrm a}$ to a Dirac structure.
It changes
the original Dirac-side resistive flow $-(f+j\mu)$ to
$-f$. The auxiliary Dirac flow is $-qe$ and vanishes
by the zero-flow factor in \eqref{eq:Rhat}, giving $qe=0$.
The relation decomposition yields the unique reconstruction
\eqref{eq:Rrecover} in both directions. Since
$\dualprod{e}{j\mu}=0$ when $qe=0$, resistive power is
preserved. The recovery map is bounded, so $\mu$
remains continuous.
\end{pf}

\begin{cor}[Resistive flow constraints]\label[cor]{cor:Rflow}
Let $\K\subseteq\Flow_{\mathrm R}$ be closed
with $\prf\Rc\subseteq\K$. Put $\Y=\Flow_{\mathrm R}/\K$,
$q:\Flow_{\mathrm R}\to\Y$, $\M=\K^\perp$, $\N=\M^\circ$,
and $j=q\adjun:\Y\adjun\to\Effort_{\mathrm R}$. Set
\[
 \widetilde\Rc=(\Rc\cap(\Flow_{\mathrm R}\times \N))
         \perp(\M\times\{0\}).
\]
On the \emph{Dirac-side} resistive flow $F_{\mathrm R}$
and auxiliary flow $\lambda$, define
\[
 B_{\mathrm f}\pvek{F_{\mathrm R}}{\lambda}
     =\pvek{-j\lambda}{qF_{\mathrm R}}.
\]
Let $\mathcal B_{\mathrm f}$ act by $B_{\mathrm f}$ on these
channels and vanish on all others. Set
\begin{align*}
 \D^{\mathrm a}&=\D\times(\Y\adjun\times\{0\}_{\Y}),\\
 \widehat\Rc&=\widetilde\Rc\times
                   (\Y\adjun\times\{0\}_{\Y}),\\
 \widehat\D&=\setdef*{(F,E+\mathcal B_{\mathrm f}F)}
                     {(F,E)\in\D^{\mathrm a}}.
\end{align*}
Then $\widehat\D$ is Dirac and $\widehat\Rc$ is
self-adjoint and nonpositive.  
Classical trajectories correspond bijectively, preserving the storage
and external port variables; the original resistive pair is uniquely
reconstructed from
\[
 f_{\mathrm R}=f,\qquad e_{\mathrm R}=e+j\lambda,
 \qquad qf=0,\quad(f,e)\in\widetilde\Rc.
\]
Here $\lambda$ denotes the auxiliary Dirac-side flow and is
uniquely determined by $\lambda=r_\K e_{\mathrm R}$.
No dynamic state is added.
\end{cor}
\begin{pf}
By \cref{lem:relaxation}, the relaxed resistor remains
self-adjoint and nonpositive. The skew-adjointness of
$B_{\mathrm f}$ ensures that the effort shear maps Dirac
structures to Dirac structures. Since $F_{\mathrm R}=-f$, the
shear adds $-j\lambda$ to its effort and $-qf$ to the
auxiliary effort. The latter vanishes, giving $qf=0$.
Decomposing the original resistive effort as $e+j\lambda$
proves unique correspondence in both directions; power
is preserved because
$\dualprod{j\lambda}{f}=\dualprod{\lambda}{qf}=0$.
\end{pf}

\noindent\emph{Example (an open-circuit constraint).}
With the resistive port convention of \eqref{eq:full-system}, an
Ohmic branch and open circuit are modeled by
\[
 \Rc=\setdef*{\left(\pvek{i}{0},\pvek{-ri}{v}\right)}
 {i,v\in\mathbb C},
 \qquad r>0.
\]
The second current vanishes while its voltage is arbitrary.
This relation is self-adjoint and nonpositive, with resistive
power $-r|i|^2$. Transferring the flow constraint with
$\K=\mathbb C\times\{0\}$ in \cref{cor:Rflow} gives
\[
 \widetilde\Rc=\setdef*{\left(\pvek{i}{i_0},\pvek{-ri}{0}\right)}
 {i,i_0\in\mathbb C}.
\]
The auxiliary Dirac condition enforces $i_0=0$, while its
multiplier reconstructs the original voltage $v=\lambda$.
Thus an ideal open-circuit constraint is moved from the
resistive relation to the power-conserving interconnection.

\begin{prop}[Storage and port power preservation]\phantom{For }
For each transfer in \cref{thm:DL,thm:LD,thm:Rtransfer,cor:Rflow},
corresponding classical trajectories have the same original state,
external input and output, stored energy, and resistive power. In
particular, both representations satisfy the same identity
\eqref{eq:power}.
\end{prop}
\begin{pf}
The external port variables are unchanged by construction. For
Dirac-to-Lagrangian, the extra storage state is zero; for
Lagrangian-to-Dirac, the extra energy term is
$\Rea\dualprod{\mu}{qx}$, which vanishes on admissible
trajectories. Resistive efforts or flows change only by elements
in the annihilator of their constrained partner. Hence
$\dualprod{e}{j\mu}=0$ or
$\dualprod{j\lambda}{f}=0$, respectively.
\end{pf}

\section{Two PDE illustrations}
Both examples involve algebraic constraints on a spatial
subregion. For simplicity, we consider one-dimensional
spatial domains. In higher dimensions, analogous transfers
can be performed using Dirac structures formulated in terms
of suitable boundary trace spaces and their dualities
(see, e.g., \citealp{Skrepek2021,BrugnoliHaineMatignon2023}
for the underlying trace theory).
In both examples, we identify $\U=\mathbb C^2$ and the respective quotient spaces $\Y$ with their anti-duals by Riesz. 
Time-dependent fields are spatial functions, e.g., $p(t)=p(t,\cdot)$; a prime denotes differentiation with respect to $\xi$.

\subsection{A string with a massless elastic segment}
Let $\Omega=(0,\ell)$ and $\Omega_0=(a,b)$, where
$0<a<b<\ell$. Put $\Omega_1=\Omega\setminus\overline\Omega_0$.
Take $\rho,T\in L^\infty(\Omega;\mathbb R)$ with
$\rho=0$ on $\Omega_0$, $\rho\geq\rho_->0$ on $\Omega_1$,
and $T\geq T_->0$ everywhere. With displacement $w$, velocity
$v=\partial_t w$, strain $\varepsilon=\partial_\xi w$, momentum
$p=\rho v$, and tension $\tau=T\varepsilon$, the model is
\begin{equation}\label{eq:wave}
 \partial_t p(t)=\partial_\xi\tau(t),\qquad
 \partial_t\varepsilon(t)=\partial_\xi v(t).
\end{equation}
The boundary input and output are endpoint forces and velocities,
\begin{equation}\label{eq:wave-ports}
 u(t)=\pvek{-\tau(t,0)}{\phantom{-}\tau(t,\ell)},\qquad
 y(t)=\pvek{v(t,0)}{v(t,\ell)}.
\end{equation}
Let $\X=\Z=L^2(\Omega;\mathbb C^2)$, with Riesz identification
of $\Z$ and $\Z\adjun$. The standard wave Stokes--Dirac structure is
\begin{equation}\label{eq:wave-D}
 \D_\partial=\setdef*{\left(
 \pvek{f_{\mathrm S}}{-y},\pvek{e_{\mathrm S}}{u}\right)}
 {v,\tau\in H^1(\Omega)},
\end{equation}
where $f_{\mathrm S}=\spvek{\tau'}{v'}$ and
$e_{\mathrm S}=\Rf_\Z\spvek{v}{\tau}$, with the boundary
ports as in \eqref{eq:wave-ports}.
Integration by parts gives isotropy, while compactly
supported tests and arbitrary boundary traces yield
self-orthogonality with respect to the symmetrized power form;
see \citep{vanderSchaftMaschke2002,LeGorrecZwartMaschke2005,
Villegas2007}. The important point is the \emph{domain condition} $v,\tau\in H^1(\Omega)$:
it enforces spatial and interface regularity, but is not itself an
algebraic constraint in $\Z$.

The storage law on $x=\spvek{p}{\varepsilon}$ is
\begin{equation}\label{eq:wave-L}
 \Lagr=\setdef*{\left(
 \spvek{\rho v}{\varepsilon},
 \Rf_\Z\spvek{v}{T\varepsilon}\right)}
 {\spvek{v}{\varepsilon}\in \Z}.
\end{equation}
The relation \eqref{eq:wave-L} has the closed state projection
\[
 \K=\setdef*{\spvek{p}{\varepsilon}\in \Z}{p|_{\Omega_0}=0}.
\]
Consequently $\Y=\Z/\K\cong L^2(\Omega_0)$ is infinite-dimensional.
Set $\alpha=\rho^{-1}$ on $\Omega_1$ and $\alpha=0$ on
$\Omega_0$, and extend $\mu\in L^2(\Omega_0)$ by zero to
$\widetilde\mu$ on $\Omega$. The relaxed storage law is
$\spvek{p}{\varepsilon}\mapsto\Rf_\Z\spvek{\alpha p}{T\varepsilon}$.
The Lagrangian-to-Dirac transfer gives the explicit new law
\[
 \widehat x=\begin{pmatrix}p\\\varepsilon\\\mu\end{pmatrix},\qquad
 \widehat e_{\mathrm S}=\pvek{\Rf_\Z\spvek{\alpha p+\widetilde\mu}{T\varepsilon}}
                   {p|_{\Omega_0}}.
\]
The auxiliary Dirac effort is zero and hence forces
$p|_{\Omega_0}=0$; thus $v=\alpha p+\widetilde\mu$.
The additional state $\mu=v|_{\Omega_0}$ is the velocity of
the massless segment, not an independent physical degree of freedom.
The endpoint ports in \eqref{eq:wave-ports} are unchanged.
More explicitly, the transformed equations are
\begin{equation}\label{eq:wave-transformed-pde}
 \begin{aligned}
  \partial_t p(t)&=\partial_\xi(T\varepsilon(t)),\\
  \partial_t\varepsilon(t)
     &=\partial_\xi(\alpha p(t)+\widetilde\mu(t)),\\
  0&=p(t)|_{\Omega_0}.
 \end{aligned}
\end{equation}
The evolution equations are unchanged in physical variables;
the last line is now a Dirac rather than a storage constraint. Crucially, the
transformed Dirac domain still requires
\[
 \alpha p(t)+\widetilde\mu(t)\in H^1(\Omega),\qquad
 T\varepsilon(t)\in H^1(\Omega).
\]
These restrictions retain the original spatial regularity and
interface coupling. A classical extension requires
$\mu\in C^1(I;L^2(\Omega_0))$, as in \cref{thm:LD}.

On the massless interval, $p(t)=0$ and \eqref{eq:wave}
gives $\partial_\xi\tau(t)=0$. Thus its tension equals a
spatially constant $F(t)$.
For strains compatible with $\varepsilon=\partial_\xi w$,
\begin{equation}\label{eq:effective-spring}
 F(t)=k_{\mathrm{eff}}[w(t,b)-w(t,a)],\qquad
 k_{\mathrm{eff}}=\left(\int_a^bT^{-1}\,\mathrm d\xi\right)^{-1}.
\end{equation}
Identity \eqref{eq:effective-spring} identifies the massless segment
as an elastic connector between two massive string sections. Its
deformation stores energy without adding inertia. The classical
power balance reads
\[
 \begin{aligned}
 \mathcal H(t)&=\frac12\int_{\Omega_1}\rho^{-1}|p(t)|^2\,\mathrm d\xi
  +\frac12\int_\Omega T|\varepsilon(t)|^2\,\mathrm d\xi,\\
 \tfrac{\mathrm d}{\mathrm dt}\mathcal H(t)
   &=\Rea\scprod{u(t)}{y(t)}_{\mathbb C^2}.
 \end{aligned}
\]

\subsection{A heat conductor with a perfectly conducting subregion}
Now let $\Omega=(0,\ell)$, with $\Omega_0=(a,b)$ and
$\Omega_1=\Omega\setminus\overline\Omega_0$. The heat capacity
$c\in L^\infty(\Omega)$ satisfies $c\geq c_->0$, while
$\kappa\in L^\infty(\Omega_1)$ satisfies
$0<\kappa_-\leq\kappa\leq\kappa_+$. For the temperature
variation $\theta$ and heat flux $\phi$, consider
\begin{equation}\label{eq:heat}
 \begin{aligned}
 c\,\partial_t\theta(t)&=-\partial_\xi\phi(t),\\
 \phi(t)&=-\kappa\partial_\xi\theta(t)
                         &&\text{on }\Omega_1,\\
 \partial_\xi\theta(t)&=0&&\text{on }\Omega_0.
 \end{aligned}
\end{equation}
The last condition describes a perfectly conducting region: the
heat flux there is not determined by a finite conductivity law.
We choose the boundary input and output
\[
 u(t)=\pvek{\theta(t,0)}{\theta(t,\ell)},\qquad
 y(t)=\pvek{\phi(t,0)}{-\phi(t,\ell)}.
\]
Let $\Z=L^2(\Omega)$, again using Riesz identification. Take
$\X=\Flow_{\mathrm R}=\Z$,
$x=c\theta$, $e_{\mathrm S}=\Rf_\Z\theta$,
$f_{\mathrm R}=\phi$, and $e_{\mathrm R}=\Rf_\Z\theta'$.
The storage law is the positive Lagrangian graph
$\Lagr=\graph(\Rf_\Z c^{-1})$.
With $g=-\phi$, the Dirac relation is given by
\begin{equation}\label{eq:heat-D}
 \D_\partial=\setdef*{\left(
 \begin{pmatrix}g'\\g\\-y\end{pmatrix},
 \begin{pmatrix}\Rf_\Z\theta\\\Rf_\Z\theta'\\u\end{pmatrix}
 \right)}{\theta,g\in H^1(\Omega)},
\end{equation}
with the boundary ports as defined above.
As for \eqref{eq:wave-D}, integration by parts shows that the Dirac structure \eqref{eq:heat-D} is isotropic,
while compactly supported tests and arbitrary endpoint traces
yield self-orthogonality; compare \citet[Thm.~6.21]{Villegas2007}.
The domain $\theta,\phi\in H^1(\Omega)$ also enforces
continuity of the physical heat flux across the internal interfaces.

For the resistive law, we identify $\Z$ with its anti-dual and regard $\K=L^2(\Omega_1)$ as a closed subspace of $\Z$ by extension by zero on $\Omega_0$.
Then
\begin{equation}\label{eq:heat-R}
 \Rc=\setdef*{(-\kappa e+h,e)}
 {e\in \K,\ h\in \K^\perp},
\end{equation}
where $\kappa e$ is extended by zero outside $\Omega_1$.
This is precisely the constitutive law in \eqref{eq:heat}.
It is self-adjoint and nonpositive. Indeed, it is the graph
of the bounded self-adjoint map
$\phi\mapsto-\kappa^{-1}P_\K\phi$ in flow-to-effort orientation.
The effort constraint $e|_{\Omega_0}=0$ is a closed, genuinely
infinite-dimensional constraint, since
$\Z/\K\cong L^2(\Omega_0)$.

The transfer in \cref{thm:Rtransfer} replaces \eqref{eq:heat-R}
by the resistive relation associated with a bounded
effort-to-flow operator,
\[
 \widetilde\Rc=\setdef*{(-\kappa P_\K e,e)}{e\in \Z}.
\]
The effort projection is now all of $\Z$.
The lossless auxiliary Dirac port has effort
$\mu\in L^2(\Omega_0)$, whose zero extension
$\widetilde\mu$ supplies the otherwise undetermined
heat flux. The transformed equations are
\begin{equation}\label{eq:heat-recover}
 \begin{aligned}
  \phi(t)&=f(t)+\widetilde\mu(t),\quad
  \mu(t)=\phi(t)|_{\Omega_0},\\
  c\,\partial_t\theta(t)
     &=-\partial_\xi(f(t)+\widetilde\mu(t)),\\
  f(t)&=-\kappa P_\K e(t),\quad e(t)=\partial_\xi\theta(t),\\
  0&=e(t)|_{\Omega_0}.
 \end{aligned}
\end{equation}
The bounded resistive law is defined on all of $\Z$,
while the last equation is imposed by the auxiliary
Dirac port. The original domain restrictions remain
\[
 \theta(t)\in H^1(\Omega),\qquad
 f(t)+\widetilde\mu(t)\in H^1(\Omega).
\]
Thus the reconstructed flux remains continuous across
$a$ and $b$; unlike \eqref{eq:wave-transformed-pde},
neither a new state nor additional temporal regularity
is required.

Since $\partial_\xi\theta(t)=0$ on $\Omega_0$, the
temperature there is spatially constant, say $\Theta(t)$. Integrating the heat balance
over $\Omega_0$ gives, for classical trajectories,
\begin{equation}\label{eq:thermal-capacity}
 C_0\dot\Theta(t)=\phi(t,a)-\phi(t,b),\qquad
 C_0=\int_a^b c(\xi)\,\mathrm d\xi.
\end{equation}
The identity \eqref{eq:thermal-capacity} shows that the ideal
conductor behaves as a lumped thermal capacity between two
distributed conductive sections. The quadratic
storage functional (not the thermodynamic internal energy) satisfies
\begin{equation}\label{eq:heat-energy}
 \begin{aligned}
 \mathcal H(t)&=\frac12\int_\Omega c|\theta(t)|^2\,\mathrm d\xi,\\
 \frac{\mathrm d}{\mathrm dt}\mathcal H(t)
 &=\Rea\scprod{u(t)}{y(t)}_{\mathbb C^2}-\int_{\Omega_1}\kappa|\partial_\xi\theta(t)|^2\,\mathrm d\xi.
 \end{aligned}
\end{equation}
Both dissipation and boundary power in \eqref{eq:heat-energy}
are preserved by the transfer.

\noindent\emph{On the opposite flow constraint.}
By \cref{cor:Rflow}, relaxing the insulating relation
$\widetilde\Rc$ recovers \eqref{eq:heat-R}.
With the original $\D_\partial$, insulation implies
$\phi|_{\Omega_0}=0$ and hence
$\partial_t\theta|_{\Omega_0}=0$.
Global $H^1$ regularity forces zero interface fluxes
and time-independent interface temperatures; avoiding
the latter requires separate temperature traces in an
adapted Dirac domain.

\section{Conclusion}
The constraint transfers between Dirac, Lagrangian, and
resistive relations preserve adjoint symmetries and power
properties without replacing unbounded operator domains
by their closures. The two PDE examples illustrate transfers
of infinite-dimensional storage and resistive constraints
to the Dirac structure, retaining spatial domain restrictions
and boundary ports.

\bibliography{Constraints}
\end{document}